\documentclass{article}

\title{\LARGE \textbf{Graph Invariants and Large Cycles - \\ a Catalog of Pure Links }}
\author{Zh.G. Nikoghosyan  \\ \\
Institute for Informatics and Automation Problems\\ National Academy of Sciences\\
P. Sevak 1, Yerevan 0014, Armenia\\ 
E-mail: zhora@ipia.sci.am}

\begin{document}

\maketitle

\begin{abstract}

Graph invariants provide a powerful analytical tool for investigation of abstract structures of graphs. They, combined in convenient relations, carry global and general information about a graph and its various substructures such as cycle structures, factors, matchings, colorings, coverings, and so on, whose discovery is the primary problem of graph theory. The major goal of this paper is to catalogue all pure relations between basic invariants of a graph and its large cycle structures, namely Hamilton, longest and dominating cycles and some their generalizations.  Basic graph invariants and pure relations allow to focus on results having no forerunners. These simplest kind of "ancestors" form a source from which nearly all possible hamiltonian results can be developed further by various additional new ideas, generalizations, extensions, restrictions and structural limitations, as well as helping researchers to make clear and simple imagination about "complicated" developmental mechanisms in the area.\\


\end{abstract}

\section{Introduction}

Graph invariants provide a powerful and maybe the single analytical tool for investigation of abstract structures of graphs. They, combined in convenient relations, carry global and general information about a graph and its particular substructures such as cycle structures, factors, matchings, colorings, coverings, and so on. The discovery of these relations is the primary problem of graph theory. Among numerous such relations there are very few and exclusive ones in forms of pure links between basic invariants of a graph and its certain substructures. Having no forerunners, these simplest kind of "ancestors" form a source from which nearly all possible results on a particular subject can be developed further by various additional new ideas, generalizations, extensions, restrictions and structural limitations.  

Hamiltonian graph theory is one of the oldest and attractive fields in graph theory concerning various path and cycle existence problems in graphs. These problems mainly are known to be NP-complete that force the graph theorists to direct efforts towards understanding the global and general relationship between various invariants of a graph and its path and cycle structure. 

The major goal of this paper is to catalogue all pure relations between basic invariants of a graph and its large cycle structures, perhaps the most important cycle structures in graphs, namely Hamilton, longest, dominating and some generalized cycles including Hamilton and dominating cycles as special cases. As in general case, these simplest kind of very few relations have no forerunners in the area actually forming a source from which nearly all possible hamiltonian results can be developed further by:  

\begin{itemize}
\item \textbf{generalized and extended graph invariants} - degree sequences (P\'{o}sa type), degree sums (Ore type, Fun-type), neighborhood unions, generalized degrees, local connectivity, and so on,
\item \textbf{extended list of path or cycle structures} - Hamilton, longest and dominating cycles, generalized $PD_\lambda$  and  $CD_\lambda$-cycles including Hamilton and dominating cycles as special cases, 2-factor, multiple Hamilton cycles, edge disjoint Hamilton cycles, powers of Hamilton cycles, $k$-ordered Hamilton cycles, arbitrary cycles, cycle systems, pancyclic-type cycle systems, cycles containing specified sets of vertices or edges,  shortest cycles, analogous path structures, and so on, 
\item \textbf{structural (descriptive) limitations} - regular, planar, bipartite, chordal and  interval graphs, graphs with forbidden subgraphs, Boolean graphs, hypercubes,  and so on,
\item\textbf{graph extensions} - hypergraphs, digraphs and orgraphs, labeled and weighted graphs, infinite graphs, random graphs, and so on.
\end{itemize}

This observation will be useful for researchers to make  a cleare imagination about developmental mechanisms in hamiltonian graph theory including the origins, current processes and future possible developments along with various research problems. 

     We refer to \cite{[8]}, \cite{[19]} and \cite{[20]} for more background and general surveys. 

All relations collected in the paper are centered around six basic well-known graph invariants, namely order $n$, size $q$, minimum degree $\delta$, connectivity $\kappa$, independence number $\alpha$ and toughness $\tau$, as well as two additional invariants $\overline{p}$ and $\overline{c}$ related to a fixed longest cycle $C$ in a graph $G$, namely the lengths of a longest path and a longest cycle in $G\backslash C$, respectively. The order $n$, size $q$ and minimum degree $\delta$ clearly are easy computable graph invariants. In \cite{[15]}, it was proved that connectivity $\kappa$ also can be determined in polynomial time. Determining the independence number $\alpha$ and toughness $\tau$  are shown in \cite{[18]} and \cite{[3]} to be $NP$-hard problems. Moreover, it was proved \cite{[3]} that for any positive rational number $t$, recognizing $t$-tough graphs (in particular 1-tough graphs) is $NP$-hard problem.

The impact of graph invariants on cycle structures gradually grows with $n$, $q$, $\delta$, $\kappa$ and $\tau$. For example, the order $n$ and size $q$ are neutral with respect to cycle structures. Meanwhile, they become more effective combined together (Theorem 1). The minimum degree $\delta$ having high frequency of occurrence in different relations is, in a sense,  a more essential invariant than the order and size, providing some dispersion of the edges in a graph. The combinations between order $n$ and minimum degree  $\delta$ become much more fruitful especially under some additional connectivity conditions. The impact of some relations on cycle structures can be strengthened under additional conditions of the type $\delta\ge\alpha\pm i$  for appropriate integer $i$. By many graph theorists, the connectivity $\kappa$  is at the heart of all path and cycle questions providing comparatively more uniform dispersion of the edges. Significant progress has been made in the area just around connectivity  $\kappa$ (section 4). An alternate connectedness measure is toughness $\tau$ - the most powerful and less investigated graph invariant introduced by Chv\'{a}tal \cite{[12]} as a means of studying the cycle structure of graphs. Chv\'{a}tal \cite{[12]} conjectured that there exists a finite constant $t_0$ such that every $t_0$-tough graph is hamiltonian. This conjecture is still open. We have omitted a number of results involving toughness $\tau$ as a parameter since they are far from being best possible. 

Large cycle structures are centered around well-known Hamilton (spanning) cycles. Other types of large cycles were introduced for different situations when the graph contains no Hamilton cycles or it is difficult to find it. Generally, a cycle $C$ in a graph $G$ is a large cycle if it dominates some certain subgraph structures  in $G$ in a sense that every such structure has a vertex in common with $C$. When $C$ dominates all vertices in $G$ then $C$ is a Hamilton cycle. When $C$ dominates all edges in $G$ then $C$ is called a dominating cycle introduced by Nash-Williams \cite{[26]}. Further, if $C$ dominates all paths in $G$ of length at least some fixed integer $\lambda$ then $C$ is a $PD_\lambda$ (path dominating)-cycle introduced by Bondy \cite{[10]}. Finally, if $C$ dominates all cycles in $G$ of length at least  $\lambda$ then $C$ is a  $CD_\lambda$ (cycle dominating)-cycle, introduced in \cite{[36]}. The existence problems of generalized $PD_\lambda$ and $CD_\lambda$-cycles are studied in \cite{[36]}.

The present catalog  includes 36 pure relations and 5 conjectures. Theorems 17, 18 and 29 present three lower bounds for the length of a longest cycle $C$ in a graph $G$ based on minimum degree $\delta$ and $\overline{p}$, $\overline{c}$ - the lengths of a longest path and a longest cycle in $G\backslash C$, respectively.  These three lower bounds for the circumference are exceptional in long cycles theory despite the fact that the idea of using $G\backslash C$ structures lies in the base of almost all proof techniques in trying to construct long cycles in graphs by the following standard procedure: choose an initial cycle $C_0$ in $G$ and try to enlarge it via path or cycle (preferably long) structures of $G\backslash C$ and connections (preferably high) between $C_0$ and $(G\backslash C_0)$-structures. Moreover, Theorem 29 provides the single lower bound for the circumference in the area involving connectivity $\kappa$ as a parameter and growing as $\kappa$ grows.

The earliest sufficient condition for a graph to be hamiltonian (Theorem 2) states that every graph with $\delta\ge n/2$ has a Hamilton cycle. Although the bound $n/2$ in Theorem 2 can not be replaced by $(n-1)/2$, it was essentially reduced to $(n+\kappa)/3$ (Theorem 19) and then to $\max\{(n+\kappa+3)/4,\alpha\}$ (Theorem 23) by incorporating new graph invariants into these bounds. We belive that these bounds can be essentially lowered further by incorporating toughness $\tau$ into the bounds. However, the bound $(n+\kappa)/3$ still remaines the lowest within graph invariants computable in polynomial time. Moreover, the second limit example in Theorems 19 shows that for each $\kappa$, the bound $(n+\kappa)/3$ can not be replaced by $(n+\kappa-1)/3$, i.e. $(n+\kappa)/3$ can not be lowered within graph invariants $n$ and $\kappa$. Observing also that among well-known basic graph invariants actually there are no other convenient ones with noticeable impact on cycle structures, one can state that very likely the bound $(n+\kappa)/3$ can not be improved in general within graph invariants computable in polynomial time. In other words, $(n+\kappa)/3$ is probably a bound between reasonable possibility and impossibility toward understanding the difficulty of $NP$-complete problems. The expression $(n+\kappa+3)/4$ in Theorem 21 presents another analogous bound concerning dominating cycles. Finally, the circumference bounds $3\delta-\kappa$ and $4\delta-\kappa-4$ in Theorems 24 and 26, respectively, are intended to be the next two bounds that can never be enlarged further within easy computable graph invariants.

The next section is devoted to necessary notation and terminology. In section 3, we discuss initial pure relations between various basic invariants of a graph and its large cycle structures having no forerunners. In fact, they are based on the order $n$, size $q$, minimum degree $\delta$ and independence number $\alpha$. Section 4 is devoted to analogous pure relations obtained from initial ones by incorporating connectivity $\kappa$ into these relations as a parameter. Finally, section 5 is devoted to a number of analogous relations under tough conditions $\tau\ge1$ and $\tau >1$.

\section{Terminology}

Throughout this article we consider only finite undirected graphs without loops or multiple edges. A good reference for any undefined terms is \cite{[11]}. We reserve $n$, $q$, $\delta$, $\kappa$ and $\alpha$ to denote the number of vertices (order), number of edges (size), minimum degree, connectivity and independence number of a graph, respectively. Each vertex and edge in a graph can be interpreted as cycles of lengths 1 and 2, respectively. A graph $G$ is hamiltonian if $G$ contains a Hamilton cycle, i.e. a cycle containing every vertex of $G$. The length $c$ of a longest cycle in a graph is called the circumference. For $C$ a longest cycle in $G$, let $\overline{p}$ and $\overline{c}$ denote the lengths of a longest path and a longest cycle in $G\backslash C$, respectively. A cycle $C^\prime$ in $G$ is a $PD_\lambda$-cycle if $|P|\le \lambda-1$ for each path $P$ in $G\backslash C^\prime$ and is a $CD_\lambda$-cycle if $|C^{\prime\prime}|\le \lambda-1$ for each cycle $C^{\prime\prime}$ in $G\backslash C^\prime$. In particular, $PD_0$-cycles and $CD_1$-cycles are well-known Hamilton cycles and $PD_1$-cycles and $CD_2$-cycles are often called dominating cycles. Let $\omega(G)$ denote the number of components of a graph $G$. A graph $G$ is $t$-tough if $|S|\ge t\omega(G\backslash S)$ for every subset $S$ of the vertex set $V(G)$ with $\omega(G\backslash S)>1$. The toughness of $G$, denoted $\tau(G)$, is the maximum value of $t$ for which $G$ is $t$-tough (taking $\tau(K_n)=\infty$ for all $n\ge 1$). 

Let $a,b,t,k$ are integers with $k\le t$. We use $H(a,b,t,k)$ to denote the graph obtained from $tK_a+\overline{K}_t$ by taking any $k$ vertices in subgraph $\overline{K}_t$ and joining each of them to all vertices of $K_b$. Let $L_\delta$ be the graph obtained from $3K_\delta +K_1$ by taking one vertex in each of three copies of $K_\delta$ and joining them each to other. For odd $n\ge 15$, construct the graph $G_n$ from $\overline{K}_{\frac{n-1}{2}}+K_\delta+K_{\frac{n+1}{2}-\delta}$, where $n/3\le\delta\le (n-5)/2$, by joining every vertex in $K_\delta$ to all other vertices and by adding a matching between all vertices in $K_{\frac{n+1}{2}-\delta}$ and $(n+1)/2-\delta$ vertices in $\overline{K}_{\frac{n-1}{2}}$. It is easily seen that $G_n$ is 1-tough but not hamiltonian. A variation of the graph $G_n$, with $K_\delta$ replaced by $\overline{K}_\delta$ and $\delta=(n-5)/2$, will be denoted by $G^*_n$.

\section{Initial relations}

We begin with a relation insuring the existence of a Hamilton cycle based on two simplest graph invariants, namely order $n$ and size $q$. \\

\noindent\textbf{Theorem 1} (Erd\"{o}s and Gallai, 1959) \cite{[16]} 

\noindent Every graph with $q\ge\frac{n^2-3n+6}{2}$ is hamiltonian.\\

Limit example: Join $K_{n-1}$ and $K_1$ by an edge.\\

The limit example shows that Theorem 1 is best possible, i.e. the size bound $(n^2-3n+6)/2$ can not be relaxed by replacing it with $(n^2-3n+5)/2$.

The earliest sufficient condition for a graph to be hamiltonian is based on the order $n$ and minimum degree  $\delta$. \\ 
 
\noindent\textbf{Theorem 2} (Dirac, 1952) \cite{[14]} 

\noindent Every graph with $\delta\ge\frac{n}{2}$ is hamiltonian.\\

Limit example: $2K_\delta+K_1$. \\

The limit example shows that the bound $n/2$ in Theorem 2 cannot be replaced by $(n-1)/2$.

A similar pure relation was established for dominating $(CD_2)$ cycles.\\

\noindent\textbf{Theorem 3} (Nash-Williams, 1971) \cite{[26]}
 
\noindent Let $G$ be a graph with $\kappa\ge 2$  and  $\delta\ge \frac{n+2}{3}$. Then each longest cycle in $G$ is a dominating cycle.\\

Limit examples: $2K_3+K_1$; $3K_{\delta-1}+K_2$; $H(1,2,4,3)$.\\

The first limit example shows that the connectivity condition $\kappa\geq 2$ in Theorem 3 can not be replaced by $\kappa\ge1$. The second example shows that the minimum degree condition $\delta\geq (n+2)/3$ can not be replaced by $\delta\geq (n+1)/2$. Finally, the third limit example shows that the conclusion "is a dominating cycle" can not be strengthened by replacing it with "is a Hamilton cycle".

Furthermore, Jung  [23] proved the third similar result concerning  $CD_3$-cycles.\\

\noindent\textbf{Theorem 4} (Jung, 1990) \cite{[23]}

\noindent Let $G$ be a graph with $\kappa\ge 3$  and  $\delta\ge \frac{n+6}{4}$. Then each longest cycle in $G$ is a $CD_3$-cycle.\\

In 2009, the author was able to find a common generalization of Theorems 2-4 by covering $CD_\lambda$-cycles for all $\lambda\ge 1$. \\

\noindent\textbf{Theorem 5} (Nikoghosyan, 2009) \cite{[36]} 

\noindent Let $G$ be a graph and $\lambda$ a positive integer. If $\kappa\ge \lambda$  and  $\delta\ge\frac{n+2}{\lambda+1}+\lambda-2$ then each longest cycle in $G$ is a  $CD_{\min\{\lambda,\delta-\lambda+1\}}$-cycle.\\

Limit examples (Theorems 4-5): $\lambda K_{\lambda+1}+K_{\lambda-1}$ $(\lambda\geq 2)$ ; $(\lambda+1)K_{\delta-\lambda+1}+K_{\lambda}$ $(\lambda\geq 1)$ ; $H(\lambda-1,\lambda,\lambda+2,\lambda+1)$ $(\lambda\geq 2)$. \\

An analogous generalization has been conjectured in \cite{[36]} for $PD_\lambda$-cycles.\\ 

\noindent\textbf{Conjecture 1} (Nikoghosyan, 2009) \cite{[36]} 

\noindent Let $G$ be a graph and $\lambda$ a positive integer. If $\kappa\ge \lambda$  and  $\delta\ge\frac{n+2}{\lambda+1}+\lambda-2$ then each longest cycle in $G$ is a  $PD_{\min\{\lambda-1,\delta-\lambda\}}$-cycle.\\

In \cite{[26]}, it was proved that the conclusion in Theorem 3 can be strengthened under additional condition $\delta\ge \alpha$.   \\

\noindent\textbf{Theorem 6} (Nash-Williams, 1971) \cite{[26]} 

\noindent Every graph with $\kappa\ge 2$ and $\delta\ge\max\{\frac{n+2}{3},\alpha\}$ is hamiltonian.\\

This theorem has been directly generalized by the following way.  \\ 

\noindent\textbf{Theorem 7} (Fraisse, 1986) \cite{[17]} 

\noindent Let $G$ be a graph and $\lambda$ a positive integer. If $\kappa\ge\lambda+1$ and $\delta\ge\max\{\frac{n+2}{\lambda+2}+\lambda-1,\alpha+\lambda-1\}$ then $G$ is hamiltonian.\\

Limit examples (Theorems 6-7): $(\lambda+1)K_{\delta-\lambda+1}+K_\lambda$ $(\delta\ge2\lambda)$; $(\lambda+2)K_{\delta-\lambda}+K_{\lambda+1}$ $(\delta\ge2\lambda+1)$; $H(\lambda,\lambda+1,\lambda+3,\lambda+2)$.\\

Now we turn to the circumference. The second earliest and simplest hamiltonian result \cite{[14]} links the circumference c and minimum degree $\delta$. \\

\noindent\textbf{Theorem 8} (Dirac, 1952) \cite{[14]} 

\noindent In every graph,  $c\ge \delta+1$.\\

Limit example: Join two copies of $K_{\delta+1}$ by an edge.\\

The same well-known paper \cite{[14]} includes the third earliest hamiltonian  relationship between minimum degree  $\delta$, circumference $c$ and Hamilton cycles. \\

\noindent\textbf{Theorem 9} (Dirac, 1952) \cite{[14]}

\noindent Let $G$ be a graph with $\kappa\ge 2$. Then $c\ge\min\{n, 2\delta\}$.\\

A similar relation has been developed for dominating cycles.\\

\noindent\textbf{Theorem 10} (Voss and Zuluaga, 1977) \cite{[38]} 

\noindent Let $G$ be a graph with  $\kappa\ge 3$. Then either $c\ge 3\delta-3$  or each longest cycle in $G$ is a dominating cycle.\\

The next common generalization covers $CD_\lambda$-cycles for all $\lambda\ge 1$ including Hamilton and dominating cycles (Theorems 9 and 10) as special cases. \\
 
\noindent\textbf{Theorem 11} (Nikoghosyan, 2009) \cite{[36]} 

\noindent Let $G$ be a graph and $\lambda$ a positive integer. If  $\kappa\ge \lambda+1$ then either $c\ge (\lambda+1)(\delta-\lambda+1)$  or each longest cycle in $G$ is a  $CD_{\min\{\lambda,\delta-\lambda\}}$-cycle. \\    

Limit examples (Theorems 9-11): $(\lambda+1)K_{\lambda+1}+K_{\lambda}$ $(\lambda\geq 1)$; $(\lambda+3)K_{\lambda-1}+K_{\lambda+2}$ $(\lambda\geq 2)$; $(\lambda+2)K_{\lambda}+K_{\lambda+1}$ $(\lambda\geq 1)$.\\

Another version of Theorem 11 was conjectured \cite{[36]} in terms of $PD_\lambda$-cycles.\\ 

\noindent\textbf{Conjecture 2} (Nikoghosyan, 2009) \cite{[36]} 

\noindent Let $G$ be a graph and $\lambda$ a positive integer. If  $\kappa\ge \lambda+1$ then either $c\ge (\lambda+1)(\delta-\lambda+1)$  or each longest cycle in $G$ is a  $PD_{\min\{\lambda-1,\delta-\lambda-1\}}$-cycle. \\

The following direct generalization includes Theorem 2  as a special case.\\
 
\noindent\textbf{Theorem 12} (Alon, 1986) \cite{[1]} 

\noindent Let $G$ be a graph and $\lambda$  a positive integer. If $\delta\ge\frac{n}{\lambda+1}$ then  $c\ge\frac{n}{\lambda}$. \\

Limit examples: $(\lambda+1)K_\lambda+K_1$; $\lambda K_{\lambda+1}$.\\

In \cite{[38]}, it was proved that the bound $\min\{n, 2\delta\}$ in Theorem 9 can be essentially enlarged under additional condition  $\delta\ge \alpha$ combined with $\kappa\ge 3$.\\

\noindent\textbf{Theorem 13} (Voss and Zuluaga, 1977) \cite{[38]}

\noindent Let $G$ be a graph with $\kappa\ge 3$ and $\delta\ge \alpha$. Then $c\ge \min\{n, 3\delta-3\}$.\\

This theorem itself has been directly generalized by the following way.  \\

\noindent\textbf{Theorem 14} (Nikoghosyan, 2009) \cite{[36]} 

\noindent Let $G$ be a graph and $\lambda$  a positive integer. If $\kappa\ge \lambda+2$  and  $\delta\ge \alpha+\lambda-1$ then  $c\ge\min\{n,(\lambda+2)(\delta-\lambda)\}$. \\

Limit examples (Theorem 13-14): $(\lambda +2)K_{\lambda +2}+K_{\lambda +1}$; $(\lambda +4)K_{\lambda}+K_{\lambda +3}$; $(\lambda +3)K_{\lambda+1}+K_{\lambda +2}$.\\

The next theorem provides a lower bound for the circumference in terms of $n, \delta, \alpha$ under the hypothesis of Theorem 3. \\

\noindent\textbf{Theorem 15} (Bauer, Morgana, Schmeichel and Veldman, 1989) \cite{[4]}

\noindent Let $G$ be a graph with $\kappa\ge 2$ and $\delta\ge \frac{n+2}{3}$. Then $c\ge \min\{n, n+\delta-\alpha\}$.\\

Limit examples: $2K_\delta+K_1$; $3K_{\delta-1}+K_2$; $K_{2\delta-2,\delta}$.\\

The first pure relation involving connectivity $\kappa$ as a parameter was developed in 1972 by showing that a remarkably simple relation $\kappa\ge\alpha$ between $\kappa$ and independence number $\alpha$ is quite sufficient for a graph to be hamiltonian. \\

\noindent\textbf{Theorem 16} (Chv\'{a}tal and Erd\"{o}s, 1972) \cite{[13]} 

\noindent Every graph with $\kappa\ge \alpha$ is hamiltonian.\\

Limit example: $K_{\delta,\delta+1}$.\\

Two lower bounds for the circumference were developed based on the fixed longest cycle $C$ in a graph $G$, minimum degree $\delta$  and some path and cycle invariants of $G\backslash C$. The first one is based on $\delta$ and $\overline{p}$  - the length of a longest path in $G\backslash C$.\\

\noindent\textbf{Theorem 17} (Nikoghosyan, 1998) \cite{[31]} 

\noindent Let $G$ be a graph and $C$ a longest cycle in $G$. Then $|C|\ge(\overline{p}+2)(\delta-\overline{p})$.  \\

The next bound is based on $\delta$  and  $\overline{c}$ - the length of a longest cycle in $G\backslash C$. \\

\noindent\textbf{Theorem 18} (Nikoghosyan, 2000) \cite{[32]} 

\noindent Let $G$ be a graph and $C$ a longest cycle in $G$. Then $|C|\ge(\overline{c}+1)(\delta-\overline{c}+1)$.\\   

Limit example (Theorems 17-18): $(\kappa+1)K_{\delta-\kappa+1}+K_\kappa$.\\

\section{Improvements via connectivity invariant}

In 1981, it was established the second hamiltonian sufficient condition involving connectivity $\kappa$ as a parameter. It can be interpreted also as the first essential improvement of Theorem 2 by incorporating connectivity $\kappa$  into the minimum degree bound without any essential limitations.\\

\noindent\textbf{Theorem 19} (Nikoghosyan, 1981) \cite{[28]} 

\noindent Every graph with $\kappa \ge 2$ and $\delta\ge \frac{n+\kappa}{3}$ is hamiltonian.\\

Limit examples: $2K_\delta+K_1$; $H(1,\delta-\kappa+1,\delta,\kappa)$ $(2\le\kappa<n/2)$.\\

A short proof of Theorem 19 was given by H\"{a}ggkvist \cite{[21]}.

An analogous relation has been established for dominating cycles which can be interpreted also as an improvement of Theorem 3.\\

\noindent\textbf{Theorem 20} (Lu, Liu and Tian, 2005) \cite{[24]} 

\noindent Let $G$ be graph with $\kappa\ge 3$  and  $\delta\ge\frac{n+2\kappa}{4}$. Then each longest cycle in $G$ is a dominating cycle.\\

Limit examples: $3K_2+K_2$; $4K_2+K_3$; $H(1,2,\kappa+1,\kappa)$.\\

The second limit example shows that for $\kappa=3$ the minimum degree bound $(n+2\kappa)/4$ in Theorem 20 can not be replaced by $(n+2\kappa-1)/4$.

Later, the bound  $(n+2\kappa)/4$ was essentially reduced to $(n+\kappa+3)/4$  without any additional limitations providing a best possible result for each $\kappa\ge 3$. \\

\noindent\textbf{Theorem 21} (Yamashita, 2008) \cite{[39]} 

\noindent Let $G$ be graph with $\kappa\ge 3$  and  $\delta\ge\frac{n+\kappa+3}{4}$. Then each longest cycle in $G$ is a dominating cycle.\\

Limit examples: $3K_{\delta-1}+K_2$; $H(2,n-3\delta+3,\delta-1,\kappa)$; $H(1,2,\kappa+1,\kappa)$.\\

In view of Theorems 19 and 21, the next generalization seems reasonable.\\

\noindent\textbf{Conjecture 3} (Yamashita, 2008) \cite{[39]} 

\noindent Let $G$ be graph and $\lambda$ an integer. If $\kappa\ge \lambda\ge2$  and  $\delta\ge\frac{n+\kappa+\lambda(\lambda-2)}{\lambda+1}$ then each longest cycle in $G$ is a $PD_{\lambda-2}$ and $CD_{\lambda-1}$-cycle.\\

The minimum degree condition $\delta\ge (n+\kappa)/3$ in Theorem 19 was essentially relaxed under additional condition  $\delta\ge\alpha$ combined with $\kappa\ge 3$. \\

\noindent\textbf{Theorem 22} (Nikoghosyan, 1985) \cite{[29]}

\noindent Every graph with $\kappa\ge 3$ and $\delta\ge \max\{\frac{n+2\kappa}{4}, \alpha\}$ is hamiltonian.\\

Limit examples: $3K_2+K_2$; $4K_2+K_3$, $H(1,2,\kappa+1,\kappa)$.\\

The second limit example shows that for $\kappa=3$ the minimum degree bound $(n+2\kappa)/4$ in Theorem 22 can not be replaced by $(n+2\kappa-1)/4$.

Later, this bound $(n+2\kappa)/4$  was reduced to $(n+\kappa+3)/4$  without any limitations providing a best possible result for each $\kappa\ge 3$.\\

\noindent\textbf{Theorem 23} (Yamashita, 2008) \cite{[39]}

\noindent Every graph with $\kappa\ge 3$ and $\delta\ge \max\{\frac{n+\kappa+3}{4}, \alpha\}$ is hamiltonian.\\

Limit examples: $3K_{\delta-1}+K_2$; $H(2,n-3\delta+3,\delta-1,\kappa)$; $H(1,2,\kappa+1,\kappa)$.\\

The first and essential improvement of Theorem 9 was achieved by incorporating connectivity  $\kappa$ into the relation without any essential limitation.\\ 

\noindent\textbf{Theorem 24} (Nikoghosyan, 1981) \cite{[28]} 

\noindent Let $G$ be a graph with  $\kappa\ge3$. Then  $c\ge\min\{n,3\delta-\kappa\}$.\\

Limit examples: $3K_{\delta-1}+K_2$; $H(1,\delta-\kappa+1,\delta,\kappa)$.\\

A simple proof of Theorem 24 was given in \cite{[25]}.

An analogous relation was developed concerning dominating cycles which can be considered also as an improvement of Theorem 10.\\

\noindent\textbf{Theorem 25} (Nikoghosyan, 2009) \cite{[37]} 

\noindent Let $G$ be a graph with  $\kappa\ge 4$. Then either $c\ge 4\delta-2\kappa$   or $G$ has a dominating cycle. \\
     
Limit examples: $4K_2+K_3$; $5K_2+K_4$; $H(1,n-2\delta,\delta,\kappa)$.\\

Theorem 25 is sharp only for $\kappa=4$ as can be seen from the second limit example.

Further, the bound $4\delta-2\kappa$  in Theorem 25 was essentially improved to $4\delta-\kappa-4$  without any limitation providing a sharp bound for each $\kappa\ge 4$.\\

\noindent\textbf{Theorem 26} (M. Nikoghosyan and Zh. Nikoghosyan, 2009) \cite{[27]} 

\noindent Let $G$ be a graph with  $\kappa\ge 4$. Then either $c\ge 4\delta-\kappa-4$   or each longest cycle in $G$ is a dominating cycle.    \\         

Limit examples: $4K_{\delta-2}+K_3$; $H(2,\delta-\kappa+1,\delta-1,\kappa)$; $H(1,2,\kappa+1,\kappa)$.

In view of Theorems 24 and 26, the following conjecture seems reasonable.\\

\noindent\textbf{Conjecture 4}  

\noindent Let $G$ be a graph and $\lambda\ge 3$ an integer. If  $\kappa\ge \lambda+1$ then either $c\ge (\lambda+1)\delta-\kappa-(\lambda+1)(\lambda-2)$   or each longest cycle in $G$ is a $PD_{\lambda-2}$ and $CD_{\lambda-1}$-cycle.    \\

The bound $3\delta-\kappa$ in Theorem 24 was enlarged to $4\delta-2\kappa$ under additional condition $\delta\ge \alpha$ combined with $\kappa\ge 4$.\\

\noindent\textbf{Theorem 27} (Nikoghosyan, 1985) \cite{[30]} 

\noindent Let $G$ be a graph with $\kappa\ge 4$  and  $\delta\ge \alpha$. Then  $c\ge\min\{n,4\delta-2\kappa\}$. \\ 

Limit examples: $4K_2+K_3$; $H(1,n-2\delta,\delta,\kappa)$; $5K_2+K_4$.\\

The bound $4\delta-2\kappa$ in Theorem 27 is sharp for $\kappa=4$ (see the third limit example).

Furthermore, the bound  $4\delta-2\kappa$  was essentially improved to $4\delta-\kappa-4$  without any additional limitations providing a best possible result for each $\kappa\ge 4$. \\

\noindent\textbf{Theorem 28} (M. Nikoghosyan and Zh. Nikoghosyan, 2009) \cite{[27]} 

\noindent Let $G$ be a graph with $\kappa\ge 4$  and  $\delta\ge \alpha$. Then  $c\ge\min\{n,4\delta-\kappa-4\}$.  \\

Limit examples: $4K_{\delta-2}+K_3$; $H(1,2,\kappa+1,\kappa)$; $H(2,n-3\delta+3,\delta-1,\kappa)$.\\

The last relation in this section is an improvement of Theorem 18 involving connectivity  $\kappa$ as a parameter combined with  $\overline{c}$ and $\delta$  such that the bound is an increasing function of $\kappa$.  \\

\noindent\textbf{Theorem 29} (Nikoghosyan, 2000) \cite{[33]} 

\noindent Let $G$ be a graph with $\kappa\ge2$ and $C$ a longest cycle in $G$. If $\overline{c}\ge \kappa$ then $|C|\ge \frac{(\overline{c}+1)\kappa}{\overline{c}+\kappa+1}(\delta+2)$. Otherwise, $|C|\ge\frac{(\overline{c}+1)\overline{c}}{2\overline{c}+1}(\delta+2)$.  \\

Limit example: $(\kappa+1)K_{\delta-\kappa+1}+K_\kappa$.\\

In view of Theorem 29, the following seems reasonable for  $PD_\lambda$-cycles.\\

\noindent\textbf{Conjecture 5} (Nikoghosyan, 2009) \cite{[36]} 

\noindent Let $G$ be a graph with $\kappa\ge2$ and $C$ a longest cycle in $G$. If $\overline{p}\ge \kappa-1$ then $|C|\ge \frac{(\overline{p}+2)\kappa}{\overline{p}+\kappa+2}(\delta+2)$. Otherwise, $|C|\ge\frac{(\overline{p}+2)\overline{p}}{2\overline{p}+2}(\delta+2)$.  \\

\section{Toughness based relations}

In \cite{[22]}, it was proved that Dirac's condition $\delta\ge n/2$ in Theorem 2 can be slightly relaxed under additional 1-tough condition.   \\
 
\noindent\textbf{Theorem 30} (Jung, 1978) \cite{[22]} 

\noindent Every graph with $n\ge 11$, $\tau\ge 1$ and $\delta\ge\frac{n-4}{2}$ is hamiltonian.\\

Limit examples: $K_{\delta, \delta+1}$; $G^*_n$.\\

This bound $(n-4)/2$ itself was lowered further to $(n-7)/2$ under stronger conditions  $n\ge 30$ and  $\tau > 1$.  \\

\noindent\textbf{Theorem 31} (Bauer, Chen and Lasser, 1991) \cite{[2]} 

\noindent Every graph with $n\ge 30$,  $\tau > 1$ and $\delta\ge \frac{n-7}{2}$ is hamiltonian. \\
 
Limit examples: Non hamiltonian graph $(n=7)$ with $\tau=1$; Petersen graph.\\

Further, it was proved that the condition $\delta\ge (n+2)/3$ in Theorem 3 can be slightly relaxed under stronger 1-tough condition instead of $\kappa\ge 2$.\\

\noindent\textbf{Theorem 32} (Bigalke and Jung, 1979) \cite{[9]}

\noindent Let $G$ be a graph with  $\tau\ge 1$ and  $\delta\ge\frac{n}{3}$. Then each longest cycle in $G$ is a dominating cycle.\\

Limit examples: $2(\kappa+1)K_2+\kappa K_1$; $L_3$; $G^*_n$.\\

Theorem 15 was improved by the same way. \\

\noindent\textbf{Theorem 33} (Bauer, Schmeichel and Veldman, 1987) \cite{[7]}

\noindent Let $G$ be a graph with $\tau\ge 1$ and $\delta\ge \frac{n}{3}$. Then $c\ge \min\{n, n+\delta-\alpha+1\}$.\\

Limit examples: $K_{\delta,\delta+1}$; $L_\delta$; $G^*_n$.\\

The next theorem is a slight improvement of Theorem 6 for 1-tough graphs.   \\

\noindent\textbf{Theorem 34} (Bigalke and Jung, 1979) \cite{[9]}

\noindent Every graph with $\tau\ge 1$ and $\delta\ge\max \{\frac{n}{3},\alpha-1\}$ is hamiltonian. \\

Limit examples: $K_{\delta,\delta+1}$ $(n\ge 3)$; $L_\delta$ $(n\ge 7)$; $K_{\delta,\delta+1}$ $(n\ge 3)$.\\

Furthermore, the bound $(n+\kappa)/3$ in Theorem 19 was slightly lowered to $(n+\kappa-2)/3$.\\

\noindent\textbf{Theorem 35} (Bauer and Schmeichel, 1991) \cite{[6]} 

\noindent Every graph with $\tau\ge 1$ and $\delta\ge \frac{n+\kappa-2}{3}$ is hamiltonian.\\

Limit examples: $K_{\delta,\delta+1}$; $L_\delta$.\\

Finally, for 1-tough graphs the bound $2\delta$ in Theorem 9 was enlarged to $2\delta+2$.\\

\noindent\textbf{Theorem 36} (Bauer and Schmeichel, 1987) \cite{[5]}

\noindent Let $G$ be a graph with $\tau\ge 1$. Then $c\ge\min\{n, 2\delta+2\}$.\\

Limit examples: $K_{\delta,\delta+1}$; $L_2$.

\end{document}